\documentclass[11pt, oneside]{amsart}

\usepackage{amsaddr}
\title{On an estimate on Götzky's domain}

\author{D\'avid 
T\'oth}
\address{HUN-REN Alfr\'ed R\'enyi Institute of Mathematics,\\
Budapest, Hungary\\and\\Department of Computer Science and Information Theory\\
Budapest University of Technology and Economics}
\email{toth.david@renyi.hu}

\thanks{
The research reported in this paper is supported by the Ministry of
Innovation and  Technology and the National Research, Development and
Innovation Office  within the Artificial Intelligence National
Laboratory of Hungary.
Project no. TKP2021-NVA-02 has been implemented with the support
provided by the Ministry of Culture and Innovation of Hungary from the
National Research, Development and Innovation Fund, financed under the
TKP2021-NVA funding scheme.
The research is also supported by the MTA–RI Lendület "Momentum" Analytic Number Theory and Representation Theory Research Group and by the NKFIH (National Research, Development and Innovation Office) grant FK 135218.
}

\usepackage[utf8]{inputenc}
\usepackage[a4paper, total={6in, 9in}]{geometry}
\usepackage{hyperref}
\usepackage{amsmath}
\usepackage{amsthm}
\usepackage{amssymb}
\usepackage{indentfirst}
\usepackage{t1enc}
\usepackage{array}
\usepackage[mathscr]{eucal}
\usepackage[utf8]{inputenc} 
\usepackage[T1]{fontenc} 
\usepackage{graphpap}
\usepackage[symbols,nogroupskip,sort=def]{glossaries-extra}
\usepackage{imakeidx}
\usepackage{xcolor}

\theoremstyle{plain}
\newtheorem{thm}{Theorem}[section]

\newtheorem{prop}[thm]{Proposition}
\newtheorem{lemma}[thm]{Lemma}

\newcommand{\N}{\mathbb{N}}
\newcommand{\R}{\mathbb{R}}

\newcommand{\Q}{\mathbb{Q}}
\newcommand{\Z}{\mathbb{Z}}

\newcommand{\OO}{\mathcal{O}}
\newcommand{\HH}{\mathbb{H}}

\newcommand{\ds}{\displaystyle}

\newcommand{\eps}{\varepsilon}

\newcommand{\abs}[1]{\left|#1\right|}

\newcommand{\ar}{\begin{array}{c}}
\newcommand{\ra}{\end{array}}
\newcommand{\arr}{\begin{array}{cc}}

\newcommand{\mtx}[4]{\left[\begin{array}{cc}#1 & #2 \\ #3 & #4\end{array}\right]}

\newcommand{\tr}{\textrm{tr}\,}

\newcommand{\UU}{\mathcal{U}}
\newcommand{\TTT}{\mathcal{T}}
\newcommand{\SSS}{\mathcal{S}}

\newenvironment{gx_extremum_pf}
{\textsl{Proof of Proposition \ref{gx_extremum}.}}
{\hfill\(\Box\)}

\newenvironment{hx_extremum_pf}
{\textsl{Proof of Proposition \ref{hx_extremum}.}}
{\hfill\(\Box\)}%

\newenvironment{choose_p_for_Hac_right_pf}
{\textsl{Proof of Proposition \ref{choose_p_for_Hac_right}.}}
{\hfill\(\Box\)}

\newenvironment{choose_p_for_Hac_left_pf}
{\textsl{Proof of Proposition \ref{choose_p_for_Hac_left}.}}
{\hfill\(\Box\)}

\newenvironment{delta_right_pf}
{\textsl{Proof of Proposition \ref{delta_right}.}}
{\hfill\(\Box\)}

\newenvironment{delta_left_pf}
{\textsl{Proof of Proposition \ref{delta_left}.}}
{\hfill\(\Box\)}

\newenvironment{increasing_right_pf}
{\textsl{Proof of Proposition \ref{increasing_right}.}}
{\hfill\(\Box\)}

\newenvironment{decreasing_left_pf}
{\textsl{Proof of Proposition \ref{decreasing_left}.}}
{\hfill\(\Box\)}

\begin{document}

\begin{abstract}
A fundamental domain $F\subset \HH^2$ for the Hilbert modular group belonging to the quadratic number field 
$\Q(\sqrt{5})$ was constructed by Götzky almost a hundred years ago. He also gave
a lower bound for the height $y_1y_2$ of the points $(z_1,z_2)=(x_1+iy_1,x_2+iy_2)\in F$.
Later Gundlach used analogous domains and estimates for other fields 
as well to give a complete list of
totally elliptic conjugacy classes in some Hilbert modular groups, while
not long ago Deutsch analysed two of these domains by numerical computations and
stated some conjectures about them. We prove one of these by giving
a sharp lower bound for the height of the points of Götzky's domain.
\end{abstract}

\maketitle

\section{Introduction}

\subsection{Hilbert modular groups}
The Hilbert modular groups are fundamental examples of discrete subgroups 
of the group $G=PSL(2,\R)^n$, where $n\geq 2$. Though our focus will be on a special
case where $n=2$ holds, we shortly recall their general definition here.
Let $\Q\leq K$ be a totally real 
finite extension of the rationals and $\OO_K$ be the ring of integers in
$K$. The corresponding Hilbert modular group is defined
as
\[
\Gamma_K:=\left\{\left(\mtx{ a^{(1)} }{ b^{(1)} }{ c^{(1)} }{ d^{(1)} }, 
\dots, \mtx{ a^{(n)} }{ b^{(n)} }{ c^{(n)} }{ d^{(n)} }\right):\,
\mtx{ a }{ b }{ c }{ d }\in PSL(2,\OO_K)\right\},
\]
where $K^{(1)},\dots,K^{(n)}$ are the different embeddings of $K$ into $\R$, and
the images of an element $a\in K$ by these embeddings  are $a^{(1)},\dots,a^{(n)}$.
Once the (ordered list of) $n$ embeddings  are fixed, any element of $\Gamma_K$ 
can and will be represented by a $2\times 2$ matrix with entries in $\OO_K^{(1)}$.

The group $G$ and hence also $\Gamma_K$ act 
on the product $\HH^n$ of $n$ copies of the complex 
upper half-plane $\HH$ coordinate-wise,
its action is described by the usual action of the coordinates on $\HH$. That is,
if $\gamma = \mtx{a}{b}{c}{d}\in PSL(2,\R)$ and $z\in \HH$, then
\begin{align*}
\gamma z = \frac{az+b}{cz+d}.
\end{align*}

\subsection{Götzky's domain}

A \emph{fundamental domain} for $\Gamma_K$  is a measurable set $F\subset \HH^n$ 
such that 
$\HH^n=\bigcup_{\gamma\in\Gamma_K}\gamma(F)$ holds and (apart from a possible exceptional
set of measure zero) no two points of $F$ are on the same $\Gamma_K$-orbit. Such a domain
is described for any Hilbert modular group in \cite{siegel}. Another fundamental domain 
was constructed
by  Götzky in \cite{gotzky} for the quadratic field $K=\Q(\sqrt{5})$. Although this 
latter construction works for any Euclidean quadratic field, in general it is only proved to   
contain a fundamental domain (the proof of
this latter statement for $\Q(\sqrt{5})$ in \cite{gotzky}
works for any quadratic Euclidean field).

The shape of such domains for quadratic fields
 have been studied 
 by Cohn \cite{cohn1965numerical, cohn1965shape},
  Deutsch \cite{deutsch2002computational, deutsch},  Jespers, Kiefer and del
Río \cite{jespers2016presentations}, and Quinn and Verjovsky \cite{quinn2020cusp}.
 More recently, 
a general
reduction algorithm was given 
by Strömberg \cite{stromberg2022reduction}
for Hilbert modular groups 
over arbitrary totally real number fields.

From now on we
restrict ourselves to totally real Euclidean 
quadratic 
extensions of $\Q$. Let $K=\Q(\sqrt{d})$ be such a 
field, where $d$ is a square-free integer greater than $1$. It will be convenient to assume
that $K$ is embedded in $\R$ and then
the ring of integers in $K$ is given by
$\OO_K=\{n+m\alpha_d:\,n,m\in\Z\}$ where 
$\alpha_d=\sqrt{d}$ if $d\equiv 2,3$ modulo $4$ and
$\alpha_d=\frac{1+\sqrt{d}}{2}$ if $d\equiv 1$ modulo $4$.
Here and in the following we always take the positive square root of a positive number.

To give Götzky's domain explicitly, we introduce
the notations
\[
S_1=\mtx{1}{1}{0}{1},\qquad 
S_{\alpha_d}=\mtx{1}{\alpha_d}{0}{1},\qquad
T=\mtx{0}{-1}{1}{0},\qquad U_d=\mtx{\eps_d}{0}{0}{\eps_d^{-1}},
\]
where $\eps_d>1$ is the generator of the unit group 
$\OO_K^\times/\{\pm1\}$. Then $\Gamma_K$ is generated by the
elements represented by the matrices above.
This is proved for $\Q(\sqrt{5})$ in \cite{gotzky} and
a slight modification of that proof gives the statement
for any totally real Euclidean quadratic field.

The coordinates of a point $z\in\HH^2$ will be written as $z_k=x_k+y_ki$ ($k=1,2$). Let us define the sets
\[
\UU_d=\{z\in\HH^2:\,\eps_d^{-2}\leq y_2/y_1<\eps_d^2\},\qquad
\TTT=\{z\in\HH^2:\,\abs{z_1z_2}\geq1\},
\]
here $\UU_d$  is clearly a fundamental domain for the subgroup generated by $U_d$, while
$\TTT$ is the fundamental domain of the
$2$ element group generated by $T$. 

Next we construct a fundamental domain for the subgroup $N_d=\left<S_1,S_{\alpha_d}\right>$ consisting of
all elements of the form $\mtx{1}{\nu}{0}{1}$ where
$\nu\in\OO_K$. The action of an element of this form on the point
$z\in\HH^2$ does not change
the values $y_1$ and $y_2$, i.e.  for any fixed $s_1,s_2>0$ 
the group $N_d$ acts on the set $H_{s_1,s_2}=\{z\in\HH^2:\,y_1=s_1,y_2=s_2\}$.
This set is homeomorphic to $\R^2$ and each $N_d$-orbit is a lattice
in it. 
From a fixed orbit we choose exactly one point $z$ such that the function $\abs{z_1z_2}$
restricted to that orbit
takes its minimal value at $z$. 
Choosing one point this way from every orbit we obtain the set $\SSS_{s_1,s_2}^d$, and then
\[
\SSS_d=\bigcup_{s_1,s_2>0}\SSS_{s_1,s_2}^d
\]
is obviously
a fundamental domain for $N_d$. 

Götzky's domain is defined as 
$\mathcal{F}_d=\UU_d\cap\TTT\cap\SSS_d$ and
- as it was mentioned before -
it is shown to be a fundamental domain in the case 
$K=\Q(\sqrt{5})$ (see \cite{gotzky}).
The shape of this domain was analyzed by numerical computations  for the
fields $\Q(\sqrt{2})$ and $\Q(\sqrt{5})$ in \cite{deutsch}, where 
some conjectures were formulated - among others -
on (strict) lower bounds for the heights $y_1y_2$ of points in $\mathcal{F}_d$.

The height plays  the same role here as the imaginary 
part of a point on the complex upper-half plane when
the action of the group $SL(2,\R)$ is considered.
In \cite{siegel} it is used to divide the fundamental domain
of a Hilbert modular group into
disjoint parts (cusp regions). While in the Euclidean
case there is only one such region, lower
bounds on the height still have some significance in such 
cases,
e.g. they are used in \cite{gundlach1964fixpunkte} 
for the computation of
totally elliptic conjugacy classes
(an element of $\Gamma_K$ is totally elliptic, if
the trace of every
coordinate of it has absolute value less than $2$), or
they can affect implied constants in other estimates,
see e.g. Lemma 1.2.5 and Lemma 2.4.1 and their proofs
in the thesis \cite{dissertation} of the author.
The aforementioned lemmata have
far-reaching consequences in  estimates
of automorphic forms, and though 
the implied constants and lower bounds on
the height are important rather from
the computational point of view, it may be desirable
to  have precise results at least in such a
classical example like  Götzky's domain.

It can be surprising at first sight that the numerical computations in \cite{deutsch} did not support actual proofs
of the conjectured estimates,
since $\mathcal{F}_d$ is a straightforward generalization 
of the standard fundamental domain for the modular group $SL(2,\Z)$ in $\HH$. But in the latter case
there is no analogue of the subgroup generated by $U_d$, and more importantly, 
the fundamental domain for the subgroup generated by the parabolic motions, i.e. 
the analogue of the set $\mathcal{S}_d$ looks substantially simpler, namely it is a strip
bounded by two hyperbolic lines. By contrast, the sets $\mathcal{S}_{s_1,s_2}^d$ look 
different for various values of $s_1$ and $s_2$ making the whole picture and the computations 
much more complicated. 

However it turns out that - at least for the field $\Q(\sqrt{5})$ - the computations can be 
simplified  around the crucial points where the minimum of the height is
taken, while crude estimates are sufficient at other places to obtain
a proof of the conjecture mentioned above. Note that parts of the following arguments also
rely 
on numerical computations. But the analytic treatment of the critical places
makes it possible to turn the numerical methods into a rigorous proof, although 
a computer was used to determine the sign of the values of some polynomial or rational 
functions at 
(finitely many)
rational places. 
\section{Estimates on Götzky's domain}

\subsection{The main result.}

 From now on we focus on the field 
 $\Q(\sqrt{5})$ examined also in
\cite{deutsch}. 
The \emph{height} of a point $z=(z_1,z_2)\in\HH^2$ is defined as the product $y_1y_2$, where 
$z_k=x_k+iy_k$ ($k=1,2$). We are going to show the following:
\begin{thm}\label{result}
    If $z\in \SSS_5\cap\TTT$, then $y_1y_2\geq\sqrt{5}/4$. The same holds
    consequently for any $z\in\mathcal{F}_5$, and in this case equality
holds if and only if $z$ is fixed by a totally elliptic element of $\Gamma_{\Q(\sqrt{5})}$
represented by a matrix of trace  $\eps_5^{\pm 1}$.
There are only finitely many points in $\mathcal{F}_5$
with this property.
\end{thm}

Our initial approach is basically the same as Götzky's in \cite{gotzky} where a
weaker bound was proved for the height in the case $d=5$. Note that a
similar argument led
to an analogous result in \cite{gundlach1964fixpunkte} for $d=2$:
\begin{lemma}\label{Gotzky_lemma}
    If $z\in \SSS_5\cap\TTT$, then $y_1y_2\geq\frac{-9+\sqrt{312}}{16}>0.54$. If $z\in\SSS_2\cap\TTT$, then $y_1y_2\geq \frac{-3+\sqrt{21}}{4}>0.3956$. 
\end{lemma}
As in the proof of the lemma, we
 are going to estimate the function
\[
f_{y_1,y_2}(x_1,x_2)=x_1^2x_2^2+x_1^2y_2^2+x_2^2y_1^2+y_1^2y_2^2=\abs{z_1z_2}^2
\]
from above  on the set
$\SSS_{s_1,s_2}^d\cap \TTT$ where $s_1,s_2>0$. 
To this end 
we will estimate on the set
\[
P_{a,d}^{s_1,s_2}=\left\{z\in\HH^2:y_1 = s_1,\,
y_2= s_2,\,
\begin{array}{c}
\ds{
	-\frac{\sqrt{d_K}}{2}\leq x_1-x_2\leq \frac{\sqrt{d_K}}{2} 
	}\\[3mm]
	-1\leq(1+a)x_1+(1-a)x_2\leq 1
\end{array}
\right\},
\]
where $a\in \R$ is a parameter and $d_K$ is the discriminant of the field $K$. 
This is a (closed) parallelogram on the plane 
$\{z\in\HH^2:\,y_1=s_1,\,y_2=s_2\}$ symmetric to the origin.
By the definition of $\SSS_{s_1,s_2}^d$ every upper bound on the former  set is 
clearly an upper bound on the latter, since if 
$z\in\SSS_{s_1,s_2}^d\cap\TTT$, then for some $\nu\in\OO_K$ we have
$(z_1+\nu,z_2+\nu')\in P_{a,d}^{s_1,s_2}$, where $\nu'$ denotes the conjugate of 
 $\nu$. 
 To simplify the notation we may
write $P_{a,d}$ or simply $P_a$ instead of $P_{a,d}^{s_1,s_2}$.

Both of the estimates listed in the lemma follow exactly the same way, estimating
the terms $x_1^2x_2^2$ and $x_1^2y_2^2+x_2^2y_1^2$ on $P_{a,d}^{s_1,s_2}$
separately. These results are the best ones that can be reached
this way, hence one has to handle all terms together to obtain a sharp bound. 
The estimation of $f_{y_1,y_2}(x_1,x_2)$ can be performed by means
of elementary calculus, but (despite the simple shape of $P_a$) 
the computations below quickly become complicated.
And though the following proof gives in principle a method that can be applied
to any quadratic Euclidean field, many small tricks that bring us through 
the numerous steps of it
fails to apply even for the field $\Q(\sqrt{2})$ (see e.g. the proof of Proposition
\ref{gx_extremum}). This does not mean that 
a similar proof cannot be performed in other cases, but such an attempt would 
probably lead to even more lengthy and tiresome calculations and the distinction
of even more cases. Nonetheless, a small part of the proof will be worked out
for a general $d$.

The key step towards the proof of Theorem \ref{result}  is
the proper 
choice 
of
the parameter $a$.
Once this problem is handled accurately 
at critical places, the other cases become easily treatable by a computer.

\subsection{Outline of the proof}

In the following we consider the numbers $y_1,y_2>0$ as parameters of the function
$f_{y_1,y_2}(x_1,x_2)$ of two variables, and the parameter $a\in(-1;1)$ will always be chosen 
according to
them, more precisely it will be a function of their ratio $c=y_2/y_1$. In the following
we fix the notations $c=y_2/y_1$ and $b=y_1y_2$ and write
\[
f_{y_1,y_1}(x_1,x_2) = x_1^2x_2^2+
(x_1^2c+x_2^2c^{-1})b +b^2.
\]

Our strategy is to choose 
the parameter $a$ so that the function
$f_{y_1,y_2}(x_1,x_2)-b^2$ takes its maximum on $P_a$ at a 
certain vertex. Let us denote this maximum 
by $g(a,b,c)$, 
we will use an estimate of the form
$g(a,b,c)\leq\lambda+\mu b$
where  $\lambda,\mu\in\R$ are suitable numbers.
Then from $z\in\SSS_d\cap\TTT$ follows
$\abs{z_1z_2}\geq1$ on the one hand, and on the other hand $(z_1+\nu,z_2+\nu')\in P_{a}$ holds for some 
$\nu\in\OO_K$, hence by the definition of $\SSS_d$ we get
\[
1\leq\abs{z_1z_2}^2\leq\abs{(z_1+\nu)(z_2+\nu')}^2=f_{y_1,y_2}(x_1+\nu,x_2+\nu')\leq g(a,b,c)+b^2
\leq\lambda+\mu b+b^2,
\]
and thus
\begin{equation}\label{quadratic_eq}
0\leq \lambda-1+ \mu b+b^2
\end{equation}
holds. 
 Since $b>0$ it is enough to obtain 
  numbers  $\lambda$, $\mu$ such  that 
 the roots of the quadratic polynomial on the right hand side of
 (\ref{quadratic_eq}) are real 
 and the smaller root is negative while the other one is at least the required bound.
That is, we require
\[
\mu^2\geq 4(\lambda-1) \qquad\textrm{and}\qquad
-\mu+\sqrt{\mu^2-4(\lambda-1)}\geq \frac{\sqrt{5}}{2}.
\]
If $\mu\geq 0$ holds, then the smaller root is automatically
negative and the second inequality above
is equivalent to
\begin{equation}\label{R_ineq}
R(\lambda,\mu):=11-16\lambda-4\sqrt{5}\mu\geq0.
\end{equation}
Note also that $R(\lambda,\mu)\geq0$ already implies $\mu^2\geq 4(\lambda-1)$ hence
it is enough to check (\ref{R_ineq}) once $\mu \geq 0$.

The summary of our plan is the following:
\begin{enumerate}
    \item[1.] We are going to choose the parameter $a\in(-1;1)$ so that
the maximum $g(a,b,c)$ of the function $f_{y_1,y_2}(x_1,x_2)$ on $P_a$ is taken at a vertex.
\item[2.] We are going bound $g(a,b,c)$ from above by $\lambda+\mu b$ where
$\mu\geq 0$ and $R(\lambda,\mu)\geq 0$.
\end{enumerate}
\subsection{Restriction of the parameter \texorpdfstring{$c$}{c}}
As a first step we show that it is enough to prove the statement 
of Theorem~\ref{result}
if $c\in[\eps_d^{-1};\eps_d]$.
To this end we consider the map
\[T_n:\HH^2\to \HH^2,\qquad 
(z_1,z_2)\mapsto(\eps_d^nx_1+i\eps_d^ny_1,(\eps_d^{n})'x_2+i\eps_d^{-n}y_2)
\] 
for any $n\in\Z$. 
(recall that $(\eps_d^n)'$ denotes the conjugate of $\eps_d^n$
in $\Q(\sqrt{d})$).
Note that $\abs{z_1z_2}^2=\abs{(T_nz)_1(T_nz)_2}^2$ holds
and $T_n$ takes the set 
$P_{a,d}^{s_1,s_2}$ to
\[
T_nP_{a,d}^{s_1,s_2}=\left\{
z\in\HH^2:\,
\begin{array}{c}
y_1 = \eps_d^ns_1\\
y_2= \eps_d^{-n}s_2
\end{array},\,
\begin{array}{c}
\ds{
-\frac{\sqrt{d_K}}{2}\leq \eps_d^{-n}x_1-(\eps_d^{-n})'x_2\leq\frac{\sqrt{d_K}}{2}
}\\[3mm]
-1\leq(1+a)\eps_d^{-n}x_1+(1-a)(\eps_d^{-n})' x_2\leq 1
\end{array}
\right\}.
\]
As before, if $z\in\HH^2$, 
$y_1= \eps_d^n s_1$ and $y_2 = \eps_d^{-n}s_2$,
then there is an integer $\nu\in\OO_K$ such that
$(z_1+\nu,z_2+\nu')\in T_nP_{a,d}^{s_1,s_2}$.
Indeed, for any $\nu\in\OO_K$ we have
\begin{equation}\label{p_eq_1}
\eps_d^{-n}\left(x_1+\nu\right)
-(\eps_d^{-n})'\left(
x_2+\nu'
\right)
=
\eps_d^{-n}x_1-(\eps_d^{-n})'x_2
+\eps_d^{-n}\nu
-(\eps_d^{-n}
\nu)'.
\end{equation}
If $\eps_d^{-n}\nu = A+B\alpha_d$ where $A,B\in\Z$, then $\eps^{-n}\nu-(\eps^{-n}\nu)' = B\sqrt{d_K}$ and
hence the expression in (\ref{p_eq_1}) can be shifted into the interval $[-\sqrt{d_K}/2;\sqrt{d_K}/2]$
by choosing $B$ properly.
Similarly
\begin{align*}
(1+a)\eps_d^{-n}
\left(x_1+\nu\right)
+&(1-a)(\eps_d^{-n})' 
\left(x_2+\nu'\right)=\\
&=
(1+a)\eps_d^{-n} x_1+(1-a)(\eps_d^{-n})' x_2
+2A+B\cdot\tr \alpha_d+aB\sqrt{d_K},
\end{align*}
so this value can be shifted into the interval $[-1;1]$ by choosing $A$ independently from $B$.

Let $z\in\SSS_d\cap\TTT$ be an arbitrary point,
then $c\in[\eps_d^{2k-1};\eps_d^{2k+1}]$ for some $k\in\Z$. 
There is a $\nu\in\OO_K$ such that $(z_1+\nu,z_2+\nu')\in T_{-k}P_{a,d}^{\eps_d^{k}y_1,\eps_d^{-k}y_2}$, and if $N(z)=N(z_1,z_2)=\abs{z_1z_2}^2$, 
then (since 
$z\in\SSS_d\cap\TTT$)  we get 
\[
1\leq \abs{z_1z_2}^2\leq\abs{(z_1+\nu)(z_2+\nu')}^2=N(z_1+\nu,z_2+\nu')=N(T_{k}(z_1+\nu,z_2+\nu')).
\]
As 
$T_{k}(z_1+\nu,z_2+\nu')\in P_{a,d}^{\eps_d^{k}y_1,\eps_d^{-k}y_2}$ and
the map $T_{k}$ does not change the value $y_1y_2$, it is enough
to estimate on this parallelogram. In other words, we can and will assume
from now on that $c\in[\eps_d^{-1};\eps_d]$.
\subsection{Proof in the neighborhoods of the endpoints.}\label{endpnt_proof_section}

From here we restrict ourselves to the case $d=5$ and the parameter $d$
will mostly be omitted from the notations 
(e.g. we always write $\eps$ instead of
$\eps_d$).  
Note that we have
$d_K=5$, $\eps=\frac{1+\sqrt{5}}{2}$ and 
$\eps^{-1}=-\eps'=\frac{\sqrt{5}-1}{2}$ in this case.

In this section we prove the theorem in the 
cases when $c$ is close enough to one of the
endpoints of the interval $[\eps^{-1};\eps]$. 
Observe first that even though the points $(z_1,z_2)=(x_1+iy_1,x_2+iy_2)$
considered in the following are not necessarily contained in $\SSS_5\cap\TTT$,
yet for
any of them holds that there are points in $\SSS_5\cap\TTT$ with the same $y_1$ and $y_2$.
This means that the statement of
Lemma \ref{Gotzky_lemma} holds for them, i.e. we always have the bound $b>0.54$.

First note that
the function $f_{y_1,y_2}(x_1,x_2)$ restricted to the set $P_{a}$ 
takes its maximum on the boundary of the parallelogram, since at every local minimum or maximum in  the 
interior of $P_a$
the partial derivatives must vanish:
\[
\ar
\partial_1f_{y_1,y_2}(x_1,x_2) = 2x_1x_2^2+2x_1y_2^2 = 0,\\[3mm]
\partial_2f_{y_1,y_2}(x_1,x_2) = 2x_1^2x_2+2x_2y_1^2 = 0.
\ra
\]
As $y_1$ and $y_2$ are positive, this implies that $x_1=x_2=0$, and
at this point  $f$ clearly takes its minimum.
Moreover, since $f_{y_1,y_2}(x_1,x_2)= f_{y_1,y_2}(-x_1,-x_2)$, it
is enough to estimate on the lines
$x_1=x_2-\frac{\sqrt{d_K}}{2}$ and $x_1=-\frac{1}{1+a}-\frac{1-a}{1+a}x_2$ between
the vertices of the parallelogram. 
Here $f$ depends on only one variable, say $x:=x_2$, 
and we also omit the constant term $y_1^2y_2^2$ for now, that is, we are looking for
the maximum of 
\[
g_{y_1,y_2}(x)=f_{y_1,y_2}\left(x-\frac{\sqrt{d_K}}{2},x\right)-y_1^2y_2^2
\]
on the interval $x\in\left[\frac{\sqrt{d_K}(1+a)-2}{4};\frac{\sqrt{d_K}(1+a)+2}{4}\right]$
and the maximum of
\[
h_{y_1,y_2}(x)=f_{y_1,y_2}\left(-\frac{1}{1+a}-\frac{1-a}{1+a}x,x\right)-y_1^2y_2^2
\]
on the interval $x\in\left[\frac{-\sqrt{d_K}(1+a)-2}{4};\frac{\sqrt{d_K}(1+a)-2}{4}\right]$.
We are going to show that for an appropriate choice of the parameter $a$
both of these functions take their maximum at an endpoint of these intervals,
and since 
$f_{y_1,y_2}(x_1,x_2)=f_{y_1,y_2}(-x_1,-x_2)$ holds, in this case it is enough to
consider the maximum of $g_{y_1,y_2}(x)$. 
The proofs of the following and the latter propositions of this section are obtained 
by means of elementary analysis of polynomial functions and (the sketches of them) 
are postponed to
Section \ref{proof_section}.
\begin{prop}\label{gx_extremum}
The function $g_{y_1,y_2}(x)$ restricted to the interval
 $\left[\frac{\sqrt{d_K}(1+a)-2}{4};\frac{\sqrt{d_K}(1+a)+2}{4}\right]$
takes its maximum at an endpoint of the interval for any  
$a\in(-1;1)$.
\end{prop}

To obtain an analogous result for  $h_{y_1,y_2}(x)$ one must be 
careful with the choice of $a$.
\begin{prop}\label{hx_extremum}
Let us define the function
\[
H_{a}(c):=\left(\frac{1-a}{1+a}\right)^2c+\frac{1}{c}-\frac{1}{(1+a)^2}.
\]
If $H_{a}(c)\geq0$, then the function $h_{y_1,y_2}(x)$ restricted to 
$\left[\frac{-\sqrt{d_K}(1+a)-2}{4};\frac{\sqrt{d_K}(1+a)-2}{4}\right]$
takes its maximum at an endpoint of the interval.
\end{prop}

Once the condition in the previous proposition is fulfilled, it is enough to examine 
the values of $g_{y_1,y_2}$ at the endpoints of the corresponding interval.
For further simplifications we
substitute
 $w=x-\sqrt{d_K}/4$, i.e. consider the
function
\begin{align*}
g(w)&:=\left(w^2-\frac{d_K}{16}\right)^2
+\left[\left(w-\frac{\sqrt{d_K}}{4}\right)^2c+
\left(w+\frac{\sqrt{d_K}}{4}\right)^2c^{-1}\right]b
\end{align*}
 at the points $w=\frac{a\sqrt{d_K}-2}{4}$ and $w=\frac{a\sqrt{d_K}+2}{4}$.
Let  $g_1(a,b,c)=g\left(\frac{a\sqrt{d_K}+2}{4}\right)$ and 
$g_2(a,b,c)=g\left(\frac{a\sqrt{d_K}-2}{4}\right)$. To decide which
value is bigger we work with their difference:
\begin{align*}
\Delta_{a,b,c} := g_1(a,b,c)-g_2(a,b,c)  
=\frac{\sqrt{5}a(5a^2-1)}{16}+
\frac{\sqrt{5}b}{2}[(a-1)c+(a+1)c^{-1}],
\end{align*}
as it can be checked by a computation (using $d_K=5$).

For each subinterval of $[\eps^{-1};\eps]$ that we consider the parameter $a$  will always
be set so
that the sign of $\Delta_{a,b,c}$ does not change on that interval. 
Once $\Delta_{a,b,c}\geq 0$ we need to estimate the value $g_1(a,b,c)$ on
that particular interval while otherwise we work with  $g_2(a,b,c)$.
Note that for a fixed $a\in(-1;1)$ and $b>0$ $\Delta_{a,b,c}$ is a decreasing function of 
$c$ on the interval $[\eps^{-1};\eps]$.

The parameter $a$ will be chosen as a function of $c$.
We choose different functions on different subintervals of $[\eps^{-1};\eps]$, 
a constant function will do on the middle intervals,
while we have to be more precautious at the endpoints where
we use
linear functions.

Let us set $a=p(c-\eps)+\frac{1}{\sqrt{d_K}}$ on the interval
$c\in[\eps-\delta;\eps]$
for some $p>0$ and $\delta>0$,
and similarly,
we set $a=p'(c-\eps^{-1})-\frac{1}{\sqrt{d_K}}$ on the interval
 $c\in[\eps^{-1};\eps^{-1}+\eta]$ for some $p'>0$ and $\eta>0$.
While a detailed analysis will be made in the former case, we 
simply choose $p'=1$ in the latter which makes the computations less tedious and fortunately works. 

Let us explain first the case of the right endpoint.
\begin{prop}\label{choose_p_for_Hac_right}
	If $c\in[1;\eps]$ and $p\in[0.24;0.66]$, then for $a=p(c-\eps)+\frac{1}{\sqrt{5}}$
	we have $a\in(-1;1)$ and $H_a(c)\geq0$.
\end{prop}

We will  
choose $p$ such that 
$\Delta_{a,b,c}$ is non-negative for any $1\leq c\leq \eps$: 
\begin{prop}\label{delta_right}
	If $c\in[1;\eps]$ and $a=p(c-\eps)+1/\sqrt{5}$ where $p=0.9/\sqrt{5}$, then $\Delta_{a,b,c}\geq0$.
\end{prop}

In summary: in a neighborhood of $\eps$ with the choice $a = p(c-\eps)+1/\sqrt{5}$ where 
$p=0.9/\sqrt{5}$
the maximal value of the function $f_{y_1,y_2}(x_1,x_2)-b^2$ on $P_a$ is $g_1(a,b,c)$.
Substituting the value of $a$ in $g_1(a,b,c)$ and using the notation $q=\sqrt{5}p$ 
we get
that $g_1(a,b,c)$ is
\[
\left(\left(\frac{q(c-\eps)+3}{4}\right)^2 -\frac{5}{16}\right)^2
 +\left[\left(\frac{q(c-\eps)+3}{4}-\frac{\sqrt{5}}{4}\right)^2c+
\left(\frac{q(c-\eps)+3}{4}+\frac{\sqrt{5}}{4}\right)^2c^{-1}\right]b.
\]
This expression can be seen as a function of $c$ with a fixed parameter $b$, let us
denote its value by $g_1(b,c)$. In the following we assume  that $b<0.56$ (otherwise the claim of 
Theorem~\ref{result} holds) and  then  this function is strictly
increasing on some interval $[\eps-\delta;\eps]$:

\begin{prop}\label{increasing_right}
If $c\in[1.48;\eps]$, $b<0.56$ and $q=0.9$, then the derivative of $g_1(b,c)$ (with respect to $c$) is positive.
\end{prop} 
We are now in the position to finish the first part of the
proof. Since $g_1(b,c)$ is strictly increasing on $[1.48;\eps]$ we simply estimate it on this interval
by the value $g_1(b,\eps)$:
\begin{align*}
	g_1(b,c)&\leq \frac{1}{16}+\left[\frac{\eps^{-4}}{4}\eps+\frac{\eps^4}{4}\eps^{-1}\right]b=
	\frac{1}{16}+\frac{b}{4}(\eps^3+\eps^{-3})=\frac{1}{16}+\frac{\sqrt{5}}{2}b.
\end{align*}
It remains to check the inequality (\ref{R_ineq}) for $\lambda = \frac{1}{16}$ and
$\mu = \frac{\sqrt{5}}{2}$. We have $R(\lambda,\mu)=0$ so (\ref{R_ineq}) holds and
the theorem is proved in the case $c\in[1.48;\eps]$. We have also proved that
equality can only hold for $c =\eps$.

Now we turn to the case when $c$ is near to the other endpoint of the interval.
As we mentioned before we choose the parameter $a = c-\eps^{-1}-\frac{1}{\sqrt{5}}$
if $c\in[\eps^{-1};\eps^{-1}+\delta]$ for some small positive $\delta$ specified later.
Then we have the following:
\begin{prop}\label{choose_p_for_Hac_left}
	If $c\in[\eps^{-1};1]$ and $a=c-\eps^{-1}-\frac{1}{\sqrt{5}}$,
	then $a\in(-1;1)$ and $H_a(c)\geq0$ hold.
\end{prop}
\begin{prop}\label{delta_left}
	If $c\in[\eps^{-1};1]$ and $a=c-\eps^{-1}-\frac{1}{\sqrt{5}}$, 
	then $\Delta_{a,b,c}\leq0$.
\end{prop}

This means that in a neighborhood of $\eps^{-1}$ with the choice of
$a=c-\eps^{-1}-1/\sqrt{5}$ the maximal value of the function $f_{y_1,y_2}(x_1,x_2)-b^2$ on $P_a$ 
is $g_2(a,b,c)$.
If we substitute the value of $a$ in $g_2(a,b,c)$ then we get that this 
maximum is
\begin{align*}
&\left(\left(\frac{\sqrt{5}(c-\eps^{-1})-3}{4}\right)^2 -\frac{5}{16}\right)^2
\\[3mm]
&\qquad\quad 
+\left[\left(\frac{\sqrt{5}(c-\eps^{-1})-3}{4}-\frac{\sqrt{5}}{4}\right)^2c+
\left(\frac{\sqrt{5}(c-\eps^{-1})-3}{4}+\frac{\sqrt{5}}{4}\right)^2c^{-1}\right]b.
\end{align*}
Let us denote this expression by $g_2(b,c)$, for a fixed $b$ it is a function of $c$.
\begin{prop}\label{decreasing_left}
	For a fixed $0<b<0.56$ the function $g_2(b,c)$ is strictly decreasing in $c$ on the interval $[\eps^{-1};0.68]$.
\end{prop}

It follows from this that 
\[
g_2(b,c)\leq g_2(b,\eps^{-1})= \frac{1}{16}+
\frac{b}{4}(\eps^4\eps^{-1}+\eps^{-4}\eps )=\frac{1}{16}+\frac{\sqrt{5}}{2}b
\]
holds for any $c\in [\eps^{-1};0.68]$.
Then we get in the same way as before that the theorem holds for
such a $c$ 
and equality can hold only if
$c=\eps^{-1}$.

\subsection{Estimates on the middle intervals}

In this section we prove the theorem in the cases when $c\in(0.68;1.48)$.
We will divide this interval into subintervals and set a fixed constant parameter $a\in(-1;1)$ 
on each of them. To ensure that the inequality 
$H_a(c)\geq 0$ holds we will use the function 
\[
\tilde{H}_a(c)=\left(\frac{1-a}{1+a}\right)^2c+\frac{2}{3}-\frac{1}{(1+a)^2}
\]
that is clearly a lower bound for $H_a(c)$ on the interval $(0.68;1.48)$.
Since $\tilde{H}_a(c)$ is increasing, it is  enough to check $\tilde{H}_a(c)\geq 0$
at the left endpoint of each subinterval.
Similarly, to estimate $g_1$ or $g_2$ 
it will be sufficient to do this at the endpoints 
once their 
derivatives 
viewed as functions in $c$ have 
constant signs on some subintervals.
These derivatives are
\[
g'_1(a,b,c)=
\left[\left(\frac{\sqrt{5}a+2}{4}-\frac{\sqrt{5}}{4}\right)^2-
\left(\frac{\sqrt{5}a+2}{4}+\frac{\sqrt{5}}{4}\right)^2c^{-2}\right]b,
\]
\[
g'_2(a,b,c)=
\left[\left(\frac{\sqrt{5}a-2}{4}-\frac{\sqrt{5}}{4}\right)^2-
\left(\frac{\sqrt{5}a-2}{4}+\frac{\sqrt{5}}{4}\right)^2c^{-2}\right]b.
\]

As a first example we consider an interval
$[1;1+\delta)$ and set  
$a=0$. As 
$\tilde{H}_0(1)=2/3$, we get that $H_0(c)\geq0$ holds if $c\in[1;1+\delta]$. 
Next we consider the derivative of $\Delta_{a,b,c}$
with respect to $c$:
\[
\frac{\partial\Delta_{a,b,c}}{\partial c}=\frac{\sqrt{5}b}{2}[(a-1)-(a+1)c^{-2}]
\]
that is negative for any $a\in(-1;1)$ and $c>0$, hence 
the function 
$\Delta_{a,b,c}$ is strictly decreasing on
$[\eps^{-1}; \eps]$ for every $a\in(-1;1)$. 
That is, to show that $\Delta_{a,b,c}\leq 0$ on an interval
it is enough to check this at the left endpoint.
This is true for $a=0$ and $c=1$ since $\Delta_{0,b,1}=0$. 

It follows that the value $g_2(0,b,c)$ is an upper bound for the function 
$f_{y_1,y_2}(x_1,x_2)-b^2$ if $c\in[1;1+\delta]$. 
The derivative of $g_2(a,b,c)$ is again increasing (as a function of $c$)
for any $a\in(-1;1)$ and positive for $a=0$ and $c=1$, hence $g_2(0,b,c)$ 
is strictly increasing on $[1;1+\delta]$ and can be estimated from above by
its value at $c_0=1+\delta$.
Now if $z\in\SSS_5\cap\TTT$ with $c\in[1;c_0]$, then
\[
1\leq\abs{z_1z_2}^2\leq g_2(0,b,c_0)+b^2,
\]
i.e. $0\leq -1+g_2(0,b,c_0)+b^2$. It is enough then if the latter
quadratic polynomial  has real roots and the smaller one is less than $1/2$
(since $b>1/2$) while the other one is bigger than $\sqrt{5}/4$.
This is true for $c_0=1.08$ so with the choice $a=0$ the theorem is proved for 
any $c\in [1;1.08)$.
Note that on this subinterval (and also on the others defined below) $b$ turns out to be
strictly bigger than $\sqrt{5}/4$. 

In the next step we increase $a$ as much as possible
so that  
$\tilde{H}_a(c_0)\geq0$,
$\Delta_{a,b,c_0}\leq 0$ and $g_2'(a,b,c_0)>0$ hold. 
For 
simplicity
we choose numbers that can easily be written down,
hence (as in the case of $c_0$ above) we round down 
to $2$ decimal places.
For the estimate of $\Delta_{a,b,c_0}$ 
we examine the sign of
$(a-1)c+(a+1)c^{-1}$. 
This value is non-positive if and only if
\[
c^2\geq(1+a)/(1-a)\Longleftrightarrow a\leq (c^2-1)/(c^2+1).
\]
If this holds, then (since $b>1/2$) we have
\[
\Delta_{a,b,c}\leq \frac{\sqrt{5}a(5a^2-1)}{16}+
\frac{\sqrt{5}}{4}[(a-1)c+(a+1)c^{-1}] =: D(a,c).
\]
We will choose $a$ such that $a<(c_0^2-1)/(c_0^2+1)$ holds and the value
$D(a,c_0)$ is non-positive.
The value $a$ that we get this way will be denoted by $a_1$.
Once $a_1$ is chosen,  we increase $c_0$ as in the first step above as much as possible to get
$c_1$ and obtain the proof of the theorem for $c\in[c_0;c_1)$. 
Continuing in the same way 
determine the values
 $a_2<a_3<\dots$ and $c_2<c_3<\dots$ until we have
$c_n\geq 1.48$ for some $n\in\N$, in which case we stop. 
We summarize this algorithm in the following steps:
\begin{enumerate}
\item[1.] Set $a_0=0$, $c_0=1.08$ and $n=1$.
\item[2.] Choose the maximal $a_{n-1}<a_n\leq\frac{c_{n-1}^2-1}{c_{n-1}^2+1}$ so that 
at most
the first two decimal digits of $a_n$ after the decimal separator are non-zero,
furthermore 
$\tilde{H}_{a_n}(c_{n-1})\geq0$, 
$D(a_n,c_{n-1})\leq 0$ and $g_2'(a_n,b,c_{n-1})>0$ hold. 
\item[3.] Choose the maximal $c_n>c_{n-1}$ such that at most
the first two decimal digits of $c_n$ after the decimal separator are non-zero and
the smaller root of the polynomial $-1+g_2(a_n,b,c_n)+b^2$ 
is less than $1/2$ while the bigger one is greater than $\sqrt{5}/4$.
\item[3.] If $c_n\geq 1.48$, then stop.
\item[4.] $n\to n+1$ and continue with step $2$.
\end{enumerate}

The algorithm above gives the following values:
\[
\begin{array}{llcllcll}
a_1=0.07,& c_1=1.15,& & a_4=0.23,&c_4=1.32,& & a_7=0.33,&c_7=1.44,\\[1mm]
a_2=0.13,&c_2=1.21,& & a_5=0.27,&c_5=1.37, & & a_8=0.34,&c_8=1.46,\\[1mm]
a_3=0.18,&c_3=1.27,& & a_6=0.3,&c_6=1.41, & &  a_9=0.36,&c_9=1.48.\\[1mm]
\end{array}
\]
This makes the proof complete if $c\in[1;\eps]$. 

Now we examine the other half of the interval
and prove the assertion on a subinterval $[c_{-1};1)$. 
As before we need 
$\tilde{H}_a(c_{-1})\geq 0$, 
but this time $\Delta_{a,b,c}\geq0$ will be required,
so the latter inequality will be checked at the right endpoint.
 We will also need the condition
\[
c^2\leq(1+a)/(1-a)\Longleftrightarrow a\geq (c^2-1)/(c^2+1)
\]
(at the right endpoints of the subintervals).
Once this is fulfilled we get
\[
\Delta_{a,b,c}\geq \frac{\sqrt{5}a(5a^2-1)}{16}+
\frac{\sqrt{5}}{4}[(a-1)c+(a+1)c^{-1}],
\]
so it is enough to show that the right hand side is non-negative at the right endpoint. 
In accordance with this we work with the function
$g_1(a,b,c)$, its derivative 
with respect to $c$
is increasing for every $a\in(-1;1)$. 
We check that this derivative is negative at $1$ (at the right endpoint) 
and so we can
estimate by $g_1(a,b,c_{-1})$ (by the value at the left endpoint).
Hence for $z\in\SSS_5\cap\TTT $ we have
\[
1\leq\abs{z_1z_2}^2\leq g_1(a,b,c_{-1})+b^2,
\]
i.e. $0\leq -1+g_1(a,b,c_{-1})+b^2$. We choose $c_{-1}$ so that
the smaller root of the quadratic polynomial on the right hand side is
smaller than 
$1/2$ while the bigger one is greater than
$\sqrt{5}/4$.

We begin with $a=0$ and looking for $c_{-1}$.
We have already seen that $\Delta_{0,b,1} = 0$. Now $g_1'(0,b,1)<0$ also holds and 
for $c_{-1}=0.92$ the other conditions are fulfilled.
Then we decrease $a$ as much as we can so that the conditions 
$a\geq (c_{-1}^2-1)/(c_{-1}^2+1)$,
$\Delta_{a,b,c_{-1}}\geq 0$ and $g_1'(a,b,c_{-1})<0$ hold (and also $\tilde{H}_a(c_{-1})\geq 0$, otherwise we could not proceed). We get the value $a_{-1} = -0.08$ and 
continue
searching 
for the next left endpoint $c_{-2}$.
We repeat these steps until $c_{-n}\leq 0.68$ holds. This way we obtain
\[
\begin{array}{llcll}
a_0=0&c_{-1}=0.92,& & a_{-5}=-0.28&c_{-6}=0.73,\\[1mm]
a_{-1}=-0.08&c_{-2}=0.86,& & a_{-6}=-0.3&c_{-7}=0.71,\\[1mm]
a_{-2}=-0.14&c_{-3}=0.82,& & a_{-7}=-0.32&c_{-8}=0.7,\\[1mm]
a_{-3}=-0.19&c_{-4}=0.78,& &a_{-8}=-0.34&c_{-9}=0.69, \\[1mm]
a_{-4}=-0.24&c_{-5}=0.75,& & a_{-9}=-0.35&c_{-10}=0.68.\\[1mm]
\end{array}
\]
Hence the assertion follows for $c\in[\eps^{-1};\eps]$ and (together 
with the postponed computations in 
Section \ref{proof_section})
the
inequality $y_1y_2\geq\sqrt{5}/4$ is proved for any $z\in \SSS_5\cap\TTT$.

\subsection{The case of equality}
It is now clear 
that $y_1y_2\geq\sqrt{5}/4$ holds for every point
$z\in\mathcal{F}_5$ since it is a subset of $\SSS_5\cap\TTT$.
In this section we 
analyse the case when equality holds in the inequality above for the points of $\mathcal{F}_5$.
By the definition of the set $\mathcal{F}_5$
we have $\eps^{-2}\leq y_2/y_1<\eps^2$ for any point $z=(z_1,z_2)$ of it and we have seen in Section \ref{endpnt_proof_section}
that equality can hold only if $y_2/y_1 = \eps^{\pm1}$. 
If $y_1y_2=\sqrt{5}/4$ and $y_2/y_1 = \eps$, then
\[
y_1 
=\frac{1}{2}\sqrt{\frac{5-\sqrt{5}}{2}}=\frac{1}{2}\sqrt{1+\eps^{-2}},\qquad
y_2 
=\frac{1}{2}\sqrt{\frac{5+\sqrt{5}}{2}}=\frac{1}{2}\sqrt{1+\eps^{2}}.
\]
Following our argument above we see that for some $\nu\in\OO_K$ the point 
$(z_1+\nu,z_2+\nu')$ is in $P_{1/\sqrt{5}}$. As before, we have
\[
1\leq \abs{z_1z_2}^2\leq\abs{(z_1+\nu)(z_2+\nu')}\leq g_1(\sqrt{5}/4,\eps)+\frac{5}{16}=\frac{1}{16}+\frac{5}{8}+\frac{5}{16}=1,
\]
and this forces these values to be equal.
That is, the point $z$ can be translated to any of 
the vertices of the parallelogram $P_{1/\sqrt{5}}$, 
e.g. to the point
\[
\left(\frac{\eps^{-2}}{2}+\frac{i}{2}\sqrt{1+\eps^{-2}},
\frac{\eps^2}{2}+\frac{i}{2}\sqrt{1+\eps^2}\right)
\]
that is the fixed point of the element represented by
\[
A=\mtx{\eps^{-1}}{1-\eps^{-1}}{-1}{1}.
\]
If $S_\nu =\mtx{1}{\nu}{0}{1}$ for any $\nu\in\OO_K$, then $z$ is the fixed point of a totally elliptic element of $\Gamma_{\Q(\sqrt{5})}$
represented by a matrix of the form
$S_{\nu}^{-1}AS_\nu$ whose trace is $\eps^{-1}+1=\eps$
and hence $z$ is an elliptic fixed point in $\mathcal{F}_5$.

One gets in the same way in the case when $y_2/y_1 = \eps^{-1}$ that
$z$ is fixed by a totally elliptic element represented by a matrix of trace $\eps^{-1}$.
Finally, the finitely many equivalence classes of elliptic fixed points
are listed in Theorem 1 (Satz 1) of \cite{gundlach1964fixpunkte} and
one checks easily that once a fixed point in $\mathcal{F}_5$ is
fixed by an element of trace $\eps^{\pm1}$ then $y_1y_2=\sqrt{5}/4$ holds.
This completes the proof of Theorem \ref{result}.

\subsection{Proofs of some propositions.}\label{proof_section}
In the following we give the sketches of the proofs of some propositions stated in 
Section \ref{endpnt_proof_section}:\\[3mm]
\indent\begin{gx_extremum_pf}
	We consider the function $g(x):=g_{y_1,y_2}(x)$
on the interval $\left[\frac{\sqrt{d_K}(1+a)-2}{4};\frac{\sqrt{d_K}(1+a)+2}{4}\right]$.
	Furthermore, to make the computation simpler we substitute 
	$w=x-\sqrt{d_K}/4$ and look for the maximum of the function
	$\tilde{g}(w) = g(w+\sqrt{d_K}/4)$ on the interval
	$\left[\frac{a\sqrt{d_K}-2}{4};\frac{a\sqrt{d_K}+2}{4}\right]$.

    This latter function is given by the formula
	\begin{align*}
	\tilde{g}(w)
	&=\left(w^2-\frac{d_K}{16}\right)^2
	+\left[\left(w-\frac{\sqrt{d_K}}{4}\right)^2c+
	\left(w+\frac{\sqrt{d_K}}{4}\right)^2c^{-1}\right]b
	\end{align*}
 where $c=y_2/y_1$ and $b=y_1y_2$, and its 
	derivative  is
	\begin{align*}
	\tilde{g}'(w)
	&=4w^3+\left(2b(c+c^{-1})-\frac{d_K}{4}\right)w+
	\frac{b\sqrt{d_K}(c^{-1}-c)}{2}.
	\end{align*}
	Since $d_K=5$, from the inequalities $b> 0.54$ 
	and $c+c^{-1}\geq2$ it follows that
	the coefficient of $w$ above is positive, 

	hence  $\tilde{g}'$ is strictly increasing on $\R$ and takes the value $0$ only once.
 (Note that e.g.
 in the case $d=2$ one cannot argue this way once $c$ is close to $1$ since
 the value $d_K=8$ is quite large.)
 So $\tilde{g}$
	has only one local extremum, and this must be a local minimum
	since $\lim_{w\to\pm\infty}\tilde{g}(w)=\infty$.
    This means that independently of 
    the choice of 
    $a$ the function 
    $\tilde{g}$ and then also $g$ take their maximum on
    the intervals above at one of the endpoints. 
\end{gx_extremum_pf}\\[3mm]
\indent\begin{hx_extremum_pf}
	We consider the function $h_{y_1,y_2}(x)=h(x)$
 on the interval 
	$\left[\frac{-\sqrt{5}(1+a)-2}{4};\frac{\sqrt{5}(1+a)-2}{4}\right]$.
	We set the notations $\alpha = (1-a)/(1+a)$ and $\beta = 1/(2(1-a))$,
 then the substitution 
	$u=x+\beta$ gives 
	\begin{align*}
	\tilde{h}(u) = h(u-\beta) = h(x)
	&=\alpha^2(u^2-\beta^2)^2
	+\left[\alpha^2(u+\beta)^2c+(u-\beta)^2c^{-1}\right]b.
	\end{align*}
	Now
	\begin{align*}
	\tilde{h}'(u)
	&=4\alpha^2u^3+[2b(\alpha^2c+c^{-1})-4\alpha^2\beta^2]u+2\beta b(\alpha^2 c-c^{-1}).
	\end{align*}
	Here the coefficient of $u^3$ is positive, and as $\alpha\beta = \frac{1}{2(1+a)}$ the
	coefficient of $u$ is
	\[
	2\left(\frac{1-a}{1+a}\right)^2bc+2bc^{-1}-\frac{1}{(1+a)^2}\geq 
	\left(\frac{1-a}{1+a}\right)^2c+\frac{1}{c}-\frac{1}{(1+a)^2}=H_{a}(c).
	\]
 Now as in the proof of Proposition \ref{gx_extremum} one can show that
	the statement is true when $H_a(c)\geq0$.
\end{hx_extremum_pf}\\[3mm]
\indent\begin{choose_p_for_Hac_right_pf}	
	We set $a=p(c-\eps)+\frac{1}{\sqrt{5}}$ for some parameter $p$. 
 Since $\eps=\frac{1+\sqrt{5}}{2}$, $a\in(-1;1)$ holds for any $0<p\leq 1$ and $c\in[1;\eps]$.
	To fulfill the condition $H_{a}(c)\geq0$ it is enough to have
	\[
	\frac{1}{c}-\frac{1}{(1+p(c-\eps)+\frac{1}{\sqrt{5}})^2}\geq0.
	\]
	By a calculation, this is equivalent to
\[
    0\leq p^2c^2+
	\left(2\left(1+\frac{1}{\sqrt{5}}\right)p-2\eps p^2-1\right)c+
	\left(p\eps-\left(1+\frac{1}{\sqrt{5}}\right)\right)^2.
\]
	The right hand side above is a quadratic polynomial in $c$ with a positive leading
 coefficient, and to fulfill this condition 
        it is sufficient if its discriminant is negative, that is if
	\begin{align*}
	0&>
    4\eps p^2-4\left(1+\frac{1}{\sqrt{5}}\right)p+1.
	\end{align*}
 Any $p$ between the roots of the latter quadratic polynomial is a good choice, in particular one can choose any value in
	the interval $[0.24;0.66]$. 
\end{choose_p_for_Hac_right_pf}\\[3mm]
\indent\begin{delta_right_pf}
We have to show that	
	\[
	\sqrt{5}a(5a^2-1)+8b[(\sqrt{5}a-\sqrt{5})c+(\sqrt{5}a+\sqrt{5})c^{-1}]\geq0.
	\]
	Let us set $q:=\sqrt{5}p$ and $t:=\eps-c$,  
 then $\sqrt{5}a =-qt+1$ and multiplying the previous inequality by $c=\eps-t$
    we obtain
	\[
	f(t):= -qt(-qt+1)(-qt+2)(\eps-t)+8b[(-qt+1-\sqrt{5})(\eps-t)^2-qt+1+\sqrt{5}]\geq0.
	\]
 We need to show that for an appropriate choice of $q$ the inequality above
 holds for any $t\in[0;\eps^{-1}]$.
	First we prove the inequality
	\[
	\varphi(t):= (-qt+1-\sqrt{5})(\eps-t)^2-qt+1+\sqrt{5}\geq0
	\]
	for any $t\in[0;\eps^{-1}]$ and some $q$. 
A calculation gives
 \begin{align*}
	\varphi(t) 
	&=-t(qt^2+2(\eps^{-1}-\eps q)t+q(\eps^2+1)-4), 
	\end{align*}
	hence it is enough to show that $qt^2+2(\eps^{-1}-\eps q)t+q(\eps^2+1)-4\leq0$. 
	The roots of this polynomial are
	\begin{align*}
 \frac{\eps q-\eps^{-1}\pm
		\sqrt{-q^2+2q+\eps^{-2}}}{q},
	\end{align*}
	so it is enough to choose a  $q>0$ such that
    the discriminant $-q^2+2q+\eps^{-2}$  is positive, 
	one of the roots above is non-positive and the other one is greater than 
	$\eps^{-1}$. 
One checks easily that these and hence $\varphi(t)\geq 0$ hold e.g. for $0.3\leq p=q/\sqrt{5} \leq 0.49$.
 
	Since $b>1/2$,
	we have for such a $q$ that
	\[
	\tilde{f}(t):=-qt(-qt+1)(-qt+2)(\eps-t)+4\varphi(t)\leq f(t),
	\]
	and it is enough to show that for a certain $q$ the following holds for 
 any $t\in(0;\eps^{-1}]$:
	\[
	F(t):=\tilde{f}(t)/t=q(t-\eps)(qt-1)(qt-2)-4(qt^2+2(\eps^{-1}-\eps q)t+q(\eps^2+1)-4)\geq 0.
	\]
 Since $F(t)$ is a cubic polynomial, this can be checked easily.
 E.g. one can check that $F(0)\geq0$ holds if
 $q= 0.9$ and that  $F'(t)\geq 0$ for any $t\in[0;\eps^{-1}]$ in this case (the roots 
 of $F'$ are greater than $1$).
\end{delta_right_pf}\\[3mm]
\indent\begin{increasing_right_pf}
We show that the function $16c^2 g'_1(b,c)$ is positive on an interval
	$[1+r;\eps]$ for some $r\geq0$. This function is a polynomial of degree $5$ in $c$
 of the form $16c^2 g'_1(b,c)=A_1(c)+bB_1(c)$, where
	\[
		A_1(c)= \left(\left(q(c-\eps)+3\right)^2-5\right)
		\left(\frac{q(c-\eps)+3}{4}\right)qc^2
	\]
	and $B_1(c)$ is given by
	\begin{align*}
		2\left(q(c-\eps)+2\eps^{-2}\right)qc^3
		+\left(q(c-\eps)+2\eps^{-2}\right)^2c^2
		+2\left(q(c-\eps)+2\eps^2\right) qc-
		\left(q(c-\eps)+2\eps^2\right)^2.
	\end{align*}
We prove the statement in two steps. First, one can show easily that $B_1(c)<0$ if $c\in[1;\eps]$,
where it can be estimated from above by
\begin{align*}
    2\left(q(c-\eps)+2\eps^{-2}\right)q\eps^3
	+\left(q(c-\eps)+2\eps^{-2}\right)^2\eps^2
	+2\left(q(c-\eps)+2\eps^2\right) q\eps-
	\left(q(c-\eps)+2\eps^2\right)^2.
\end{align*}
An elementary analysis of the latter quadratic polynomial of $c$ shows that this upper bound
is still negative on $[1;\eps]$, we omit the detailed computation.

By our assumption on $b$ we have
	\[
	16c^2g'_1(b,c)=A_1(c)+bB_1(c)> A_1(c)+0.56 B_1(c).
	\]
 for $c\in [1;\eps]$,
	hence it is enough to show that the right hand side above is positive for any $c\in[1.48;\eps]$.
	To avoid the work with complicated algebraic expressions we check this
    via (finitely many) substitutions. (This method will also be applied  in the 
    subsequent proofs.)
	Namely, we show that the polynomial $F_1(c)=A_1(c)+0.56 B_1(c)+8.001$
	has $5$ roots, and therefore if $x_0$ is the biggest root, then $F_1$ 
	must be strictly increasing on the interval $[x_0;\infty)$. 
	We do all this by giving pairs $c_1,c_2$ of real numbers such that
	$c_1<c_2$ and the signs of $F_1(c_1)$ and $F_1(c_2)$ are different.
	One checks easily (e.g. by a computer)
	that 
	\[
	F_1(-14)<0,\quad F_1(-13)>0,\quad F_1(-0.1)<0,\quad F_1(0)>0,\quad F_1(0.1)<0,
	\qquad F_1(0.6)>0.
	\]
	Thus the function $A_1(c)+0.56B_1(c)$ is strictly  increasing for $c\geq0.6$. On the other hand, 
	for $c=1.48$ its value is positive  and hence the same is true for $c\geq1.48$.
\end{increasing_right_pf}\\[3mm]
\indent\begin{choose_p_for_Hac_left_pf}
One checks easily that $a\in(-1;1)$ holds.	
We have to show that
	\[
	\left(\frac{1-a}{1+a}\right)^2c+\frac{1}{c}-\frac{1}{(1+a)^2}\geq0
	\]
 holds if $c\in[\eps^{-1};1]$.
	Multiplying by $(1+a)^2c$ we get 
	\begin{align*}
	(1-a)^2c^2+(1+a)^2-c
	&\geq \left(\eps^{-1}+\frac{1}{\sqrt{5}}\right)^2c^2-c+
	\left(1-\frac{1}{\sqrt{5}}\right)^2.
	\end{align*}
	The discriminant of this latter quadratic polynomial is negative so its value
 is positive for
 any $c\in\R$ and $H_a(c)\geq0$ follows.
\end{choose_p_for_Hac_left_pf}\\[3mm]
\indent\begin{delta_left_pf}
	We  show that $16\Delta_{a,b,c}\leq 0$, i.e.
	\begin{align*}
	(\sqrt{5}(c-\eps^{-1})-1)\sqrt{5}&(c-\eps^{-1})(\sqrt{5}q(c-\eps^{-1})-2)+\\[3mm]
	&+8b[(\sqrt{5}(c-\eps^{-1})-1-\sqrt{5})c+(\sqrt{5}(c-\eps^{-1})-1+\sqrt{5})c^{-1}]\leq 0.
	\end{align*}
	Multiplying by $c$ and substituting 
	$t=c-\eps^{-1}$ we get
	\begin{align*}
	f(t)= \sqrt{5}t(\sqrt{5}t-1)&(\sqrt{5}t-2)(t+\eps^{-1})
 +8b[(\sqrt{5}t-1-\sqrt{5})(t+\eps^{-1})^2+\sqrt{5}t-1+\sqrt{5}],
	\end{align*}
	hence it must be  shown that $f(t)\leq 0$ if $t\in[0;1-\eps^{-1}]=[0;\eps^{-2}]$.
    First we check that
	\begin{equation}\label{phi_non_positive}
	\varphi(t) := (\sqrt{5}t-1-\sqrt{5})(t+\eps^{-1})^2+\sqrt{5}t-1+\sqrt{5}\leq0
	\end{equation}
	if $t\in[0;\eps^{-2}]$. 
    A computation shows that $\varphi(t)=t\tilde{\varphi}(t)$ where
 \[\tilde{\varphi}(t) = \sqrt{5}t^2-2\eps^{-3}t +\sqrt{5}(\eps^{-2}+1)-4.
 \]
 One checks that $\tilde{\varphi}(0)<0$ and
	$\tilde{\varphi}(\eps^{-2})= 
        -3+\sqrt{5}<0$,
	so $\tilde{\varphi}(t)$ is negative for any $t\in[0;\eps^{-2}]$ and 
	hence (\ref{phi_non_positive}) is proved.

	As $b>0.5$ we have
$	f(t)\leq \sqrt{5}t(\sqrt{5}t-1)(\sqrt{5}t-2)(t+\eps^{-1})+
	4\varphi(t)=: \tilde{f}(t)$.
	Since $\tilde{f}(0) = 0$, it is enough to show that $\tilde{f}'$ is negative on the interval
	$[0;\eps^{-2}]$. 
	A computation gives 
 $\tilde{f}'(t) = 20\sqrt{5}t^3+3\sqrt{5}\eps^{-2}t^2+(47-27\sqrt{5})t+9\sqrt{5}-21$.
Now the inequalities
	\[
 \tilde{f}'(-0.5)>0,\qquad
 \tilde{f}'(0)<0,\qquad
 \tilde{f}'(\eps^{-2})<0
	\]
 hold and
	 the assertion follows (because $\tilde{f}'$ is a polynomial function of degree $3$ with positive leading
	coefficient).
\end{delta_left_pf}\\[3mm]
\indent\begin{decreasing_left_pf}
	It is enough to see that 
 $16c^2 g'_2(b,c)$ is negative on the interval $[\eps^{-1}; 0.68]$.
	Similarly as earlier we write $16c^2 g'_2(b,c)=A_2(c)+bB_2(c)$ where
	\begin{align*}
	A_2(c)&= \left(\left(\sqrt{5}(c-\eps^{-1})-3\right)^2-5\right)
	\left(\frac{\sqrt{5}(c-\eps^{-1})-3}{4}\right)\sqrt{5}c^2,
 \end{align*}
 and 
 \begin{align*}
	B_2(c)&=
 2\left(\sqrt{5}(c-\eps^{-1})-2\eps^2\right)\sqrt{5}c^3
	+\left(\sqrt{5}(c-\eps^{-1})-2\eps^2\right)^2c^2
\\[3mm]
	&\quad
 +2\left(\sqrt{5}(c-\eps^{-1})-2\eps^{-2}\right) \sqrt{5}c-
	\left(\sqrt{5}(c-\eps^{-1})-2\eps^{-2}\right)^2.
	\end{align*}
Note that
 $c-\eps^{-1}\leq 0.68-\eps^{-1}<0.062$ and then
	\[
	\sqrt{5}(c-\eps^{-1})-2\eps^2<0,\qquad \sqrt{5}(c-\eps^{-1})-2\eps^{-2}<0, 
	\]
	therefore
	\begin{align*}
	B_2(c)&\geq 2\left(\sqrt{5}(c-\eps^{-1})-2\eps^2\right)\sqrt{5}\cdot 0.68^3+
	\left(\sqrt{5}(c-\eps^{-1})-2\eps^2\right)^2\eps^{-2}\\[3mm]
	&\quad+2\left(\sqrt{5}(c-\eps^{-1})-2\eps^{-2}\right) \sqrt{5}\cdot 0.68-
	\left(\sqrt{5}(c-\eps^{-1})-2\eps^{-2}\right)^2.
	\end{align*}
	This lower bound is a quadratic polynomial of  $c$ and one checks that
    it is positive on the interval $[\eps^{-1};0.68]$ and hence so is $B_2(c)$.
 We have then the upper bound
	\[
	16c^2g'_2(b,c)<A_2(c)+0.56B_2(c)
	\]
and to see that the right hand side above is negative we consider the function
$F_2(c) = A_2(c)+0.56B_2(c)+2.58$. This is a polynomial of degree $5$
	with positive leading coefficient and we have that
\[
\ar
	F_2(-0.1)<0,\quad F_2(0)>0,\quad F_2(0.2)<0,\quad F_2(0.4)>0,\\[3mm] F_2(1.2)>0,\quad 
	F_2(1.4)<0,\quad F_2(3)>0.
\ra
\]
    This implies that $F_2$ has a root $x_1$ in $[0.2;0.4]$ and another one in $[1.2;1.4]$ denoted by $x_2$.
	Furthermore, 
	$F_2$ is positive on $(x_1;x_2)$ where it has exactly one local maximum taken at
	a point $x_m$, hence $F_2$ is increasing on 
	$[x_1;x_m]$ while it is decreasing on $[x_m;x_2]$. 
	Since $F_2(0.4)<F_2(0.7)<F_2(0.8)$ we get that $x_m>0.7$ and hence $F_2$
	is increasing on the interval
	$[0.4;0.7]$ and so is
	$ A_2(c)+0.56B_2(c)$.
	Moreover, $A_2(0.7)+0.56B(0.7)<0$, therefore $g_2'(b,c)<0$ on the interval $[\eps^{-1};0.68]$.
\end{decreasing_left_pf}

\bibliographystyle{abbrv}
\bibliography{references}

\end{document}